\magnification=1200
\overfullrule=0pt
\centerline {\bf A multiplicity result for nonlocal problems}
\centerline {\bf involving
nonlinearities with bounded primitive}\par
\bigskip
\bigskip
\centerline {BIAGIO RICCERI}\par
\bigskip
\bigskip
{\bf Abstract.} The aim of this paper is to provide the first application
of Theorem 3 of [2] in a case where the dependence of the underlying
equation from the real parameter is not of affine type. The simplest
particular case of our result reads as follows:\par
\smallskip
 Let $f:{\bf R}\to {\bf R}$ be a non-zero continuous function such that
$$\sup_{\xi\in {\bf R}}|F(\xi)|<+\infty$$
where $F(\xi)=\int_0^{\xi}f(s)ds$.
Moreover, let
 $k:[0,+\infty[\to {\bf R}$ and $h:]-\hbox {\rm osc}_{\bf R}F,
\hbox {\rm osc}_{\bf R}F[\to{\bf R}$ be two continuous and non-decreasing
functions,
with $k(t)>0$ for all $t>0$ and $h^{-1}(0)=\{0\}$.\par
Then, for each $\mu$ large enough,
there exist an open interval $A\subseteq ]\inf_{\bf R}F,\sup_{\bf R}F[$
and a number $\rho>0$ such that, for every $\lambda\in A$,
the problem
$$\cases {-k\left ( \int_0^1|u'(t)|^2dt\right )u''=
\mu h\left ( \int_0^1F(u(t))dt-\lambda\right ) f(u) &
in $[0,1]$\cr & \cr u(0)=u(1)=0\cr}$$
has at least three solutions whose
norms in $H^1_0(0,1)$ are less than $\rho$.\par

\bigskip
\bigskip
\bigskip
\bigskip
In [2], we established the following result:\par
\medskip
THEOREM A. - {\it Let $X$ be a separable and reflexive real
Banach space, $I\subseteq {\bf R}$ an interval, and
$\Psi:X\times I\to {\bf R}$ a continuous
function satisfying the following
conditions:\par
\noindent
$(a_1)$\hskip 5pt for each $x\in X$, the function
$\Psi(x,\cdot)$ is concave;\par
\noindent
$(a_2)$\hskip 5pt for each $\lambda\in I$, the function
$\Psi(\cdot,\lambda)$ is $C^1$, sequentially weakly lower semicontinuous,
coercive, and satisfies the Palais-Smale condition;\par
\noindent
$(a_3)$\hskip 5pt there exists a continuous concave
function $h:I\to {\bf R}$ such that
$$\sup_{\lambda\in I}\inf_{x\in X}(\Psi(x,\lambda)+h(\lambda))<
\inf_{x\in X}\sup_{\lambda\in I}(\Psi(x,\lambda)+h(\lambda))\ .$$
Then, there exist an open interval $A\subseteq I$
and a positive real number $\rho$, such that, for each $\lambda
\in J$, the equation
$$\Psi'_{x}(x,\lambda)=0$$
has at least three solutions in $X$ whose norms are less than $\rho$.}
\medskip
A consequence of Theorem A is as follows:\par
\medskip
THEOREM B. - {\it Let $X$ be a separable and reflexive real Banach
space; $\Phi:X\to {\bf R}$ a continuously
 G\^ateaux differentiable
and sequentially weakly lower semicontinuous
 functional whose G\^ateaux derivative
admits a continuous inverse on $X^*$; $\Psi:X\to {\bf R}$ a
continuously G\^ateaux differentiable functional whose G\^ateaux
derivative is compact; $I\subseteq {\bf R}$ an interval.
 Assume that
$$\lim_{\|x\|\to +\infty}(\Phi(x)+\lambda\Psi(x))=+\infty
$$
for all $\lambda\in I$, and that there exists a
continuous concave function $h:I\to {\bf R}$ such
that
$$\sup_{\lambda\in I}\inf_{x\in X}(\Phi(x)+
\lambda\Psi(x)+h(\lambda))<\inf_{x\in X}
\sup_{\lambda\in I}(\Phi(x)+\lambda\Psi(x)+
h(\lambda))\ .$$
Then, there exist an open interval $A\subseteq I$ and a
positive real number $\rho$ such
that, for each $\lambda\in A$,
the equation
$$\Phi'(x)+\lambda\Psi'(x)=0$$
has at least three solutions in $X$ whose norms are less than $\rho$.}\par
\medskip
In appraising the literature, it is quite surprising to realize that,
while Theorem B
has been proved itself to be one of the most frequently
used abstract multiplicity results in the last decade, it seems that there is
no article where
Theorem A has been applied to some $\Psi$ which does not depend on
$\lambda$ in an affine way. For an up-dated bibliographical account related
to Theorem B, we refer to [3]. \par
\smallskip
The aim of this paper is to offer a first contribution to fill this gap.
\par
\smallskip
To state our results, let us fix some notation. \par
\smallskip
For a generic function $\psi:X\to {\bf R}$, we denote
by $\hbox {osc}_X \psi$ the (possibly infinite) number
$\sup_X\psi-\inf_X\psi$.\par
\smallskip
In the sequel, $\Omega\subset {\bf R}^n$ is a bounded domain with smooth
boundary. We consider the space
$H^1_0(\Omega)$ equipped with the norm
$$\|u\|=\left ( \int_{\Omega}|\nabla u(x)|^2dx\right )^ {1\over 2}\ .$$
If $I\subseteq {\bf R}$ is an interval, with
$0\in I$, and $g:\Omega\times I\to {\bf R}$
is a function such that $g(x,\cdot)$ is continuous in $I$ for all $x\in
\Omega$, we set
$$G(x,\xi)=\int_0^{\xi}g(x,t)dt$$
for all $(x,\xi)\in \Omega\times I$.\par
\smallskip
When $n\geq 2$,  we denote by ${\cal A}$ the class of all
Carath\'eodory functions $f:\Omega\times {\bf R}\to {\bf R}$ such that
$$\sup_{(x,\xi)\in \Omega\times {\bf R}}{{|f(x,\xi)|}\over
{1+|\xi|^q}}<+\infty\ ,$$
for some $q$ with  $0<q< {{n+2}\over {n-2}}$ if $n\geq 3$ and $0<q<+\infty$
if $n=2$.
When $n=1$, we denote by ${\cal A}$  the class
of all Carath\'eodory functions $f:\Omega\times {\bf R}\to {\bf R}$ such
that, for each $r>0$, the function $x\to \sup_{|t|\leq r}|f(x,t)|$ belongs
to $L^{1}(\Omega)$.\par
\smallskip
If $f\in {\cal A}$, for each $u\in H^1_0(\Omega)$, we set
$$J_f(u)=\int_{\Omega}F(x,u(x))dx\ .$$
The functional $J_f$ is $C^1$ and its derivative is compact.
Moreover, we set
$$\alpha_f= \inf_{H^1_0(\Omega)} J_f\ ,$$
$$\beta_f=\sup_{H^1_0(\Omega)} J_f$$
and
$$\omega_f=\beta_f-\alpha_f\ .$$
Clearly, when $f$ does not depend on $x$, we have
$$\alpha_f= \hbox {\rm meas}(\Omega)\inf_{{\bf R}} F$$
and
$$\beta_f=\hbox {\rm meas}(\Omega)
\sup_{{\bf R}} F \ .$$
Our main result reads as follows:\par
\medskip
THEOREM 1. - {\it Let $f, g\in {\cal A}$ be such that
$$\sup_{(x,\xi)\in \Omega\times {\bf R}}\max
\{ |F(x,\xi)|, G(x,\xi)\} <+\infty$$
and
$$\sup_{u\in H^1_0(\Omega)}\left | \int_{\Omega}F(x,u(x))dx\right | >0\ .$$
Then, for every pair of continuous and non-decreasing
functions $k:[0,+\infty[\to {\bf R}$ and $h:]-\omega_f,\omega_f[\to{\bf R}$, 
with $k(t)>0$ for all $t>0$ and $h^{-1}(0)=\{0\}$, for which
the number
$$\theta^*=\inf\left \{ {{{{1}\over {2}}K\left ( \int_{\Omega}|\nabla u(x)|^2dx\right )
-\int_{\Omega}G(x,u(x))dx}
\over {H\left ( \int_{\Omega}F(x,u(x))dx\right ) }}: u\in H^1_0(\Omega),
\int_{\Omega}F(x,u(x))dx\neq 0 \right \}
$$
is non-negative, and for every $\mu>\theta^*$,
there exist an open interval $A\subseteq ]\alpha_f,\beta_f[$
and a number $\rho>0$ such that, for every $\lambda\in A$,
the problem
$$\cases {-k\left ( \int_{\Omega}|\nabla u(x)|^2dx\right )\Delta u=
\mu h\left ( \int_{\Omega}F(x,u(x))dx-\lambda\right ) f(x,u)+g(x,u) &
in $\Omega$\cr & \cr u=0 & on $\partial\Omega$\cr}$$
has at least three weak solutions whose
norms in $H^1_0(\Omega)$ are less than $\rho$.}\par
\medskip
Clearly, a weak solution of the above problem
problem is any $u\in H^1_0(\Omega)$ such
that
$$k\left ( \int_{\Omega}|\nabla u(x)|^2dx\right ) \int_{\Omega}
\nabla u(x)\nabla v(x)dx=$$
$$=\mu h\left ( \int_{\Omega}F(x,u(x))dx-\lambda\right )
\int_{\Omega}f(x,u(x))v(x)dx+\int_{\Omega}g(x,u(x))v(x)dx$$
for all $v\in H^1_0(\Omega)$.\par
So, the weak solutions of the problem are exactly the critical points
in $H^1_0(\Omega)$ of the functional
$$u\to {{1}\over {2}}K(\|u\|^2)-\int_{\Omega}G(x,u(x))dx-
\mu H\left ( \int_{\Omega}F(x,u(x))dx-\lambda\right )\ .$$
\smallskip
The problem that we are considering is a nonlocal one. We refer to the
very recent paper [1] for a relevant discussion and an up-dated
bibliography as well.\par
\smallskip
  From what we said above, it is clear that our proof of Theorem 1 is based
on the use of Theorem A. This is made possible by the following
proposition:\par
\medskip
PROPOSITION 1. - {\it Let $X$ be a non-empty set and let $\gamma:
X\to [0,+\infty[$,
$J:X\to {\bf R}$ be two functions such that $\gamma(x_0)=J(x_0)=0$
for some $x_0\in X$. Moreover, assume that $J$ takes at
least four values. Finally, let $\varphi:
]-\hbox {\rm osc}_XJ,\hbox {\rm osc}_XJ[\to [0,+\infty[$ be a
continuous function such that
$$\varphi^{-1}(0)=\{0\} \eqno{(1)}$$
and
$$\min\left \{ \liminf_{t\to (-\hbox {\rm osc}_XJ)^+}\varphi(t),
\liminf_{t\to (\hbox {\rm osc}_XJ)^-}\varphi(t)\right \} >0\ .\eqno{(2)}$$
Put
$$\theta=\inf_{x\in J^{-1}(]\inf_X J, \sup_X J[\setminus \{0\})}
{{\gamma(x)}\over {\varphi(J(x))}}\ .$$
Then, for each $\mu>\theta$, we have
$$\sup_{\lambda\in ]\inf_X J, \sup_X J[}\inf_{x\in X}
(\gamma(x)-\mu\varphi(J(x)-\lambda))<
\inf_{x\in X}\sup_{\lambda\in ]\inf_X J, \sup_X J[}
(\gamma(x)-\mu\varphi(J(x)-\lambda))\ .$$}\par
\smallskip
PROOF. First, we make some remarks on the definition of $\theta$. Since
$J$ takes at least four values, the set
$J^{-1}(]\inf_X J, \sup_X J[\setminus \{0\}$ is non-empty. So, if
$x\in J^{-1}(]\inf_X J, \sup_X J[\setminus \{0\})$, we have
$J(x)\in ]-\hbox {\rm osc}_XJ,\hbox {\rm osc}_XJ[\setminus
\{0\}$ (recall that $\inf_X J\leq 0\leq\sup_X J$), and so $\varphi(J(x))>0$.
Hence, $\theta$ is a well-defined non-negative real number. Now, fix
$\mu>\theta$. Since $\varphi$ is continuous, we have
$$\inf_{\lambda\in ]\inf_XJ,\sup_XJ[}\varphi(J(x)-\lambda)=0$$
for all $x\in X$.
Hence
$$\inf_{x\in X}\sup_{\lambda\in ]\inf_X J, \sup_X J[}
(\gamma(x)-\mu\varphi(J(x)-\lambda))=\inf_{x\in X}
\left ( \gamma(x)-\mu\inf_{\lambda\in ]\inf_XJ,\sup_XJ[}
\varphi(J(x)-\lambda)\right )$$
$$ =\inf_X\gamma=0\ .\eqno{(3)}$$
Now, since $\mu>\theta$, there is $x_1\in X$ such that
$$\gamma(x_1)-\mu\varphi(J(x_1))<0\ .$$
So, again by the continuity of $\varphi$, for $\epsilon, \delta>0$
small enough, we have
$$\gamma(x_1)-\mu\varphi(J(x_1)-\lambda)<-\epsilon \eqno{(4)}$$
for all $\lambda\in [-\delta,\delta]$. On the other hand, $(1)$ 
and $(2)$ imply that
$$\nu:=\inf_{\lambda\in ]\inf_XJ,\sup_XJ[\setminus
[-\delta,\delta]}\varphi(-\lambda)>0\ .\eqno{(5)}$$
  From $(4)$ and $(5)$, recalling that $\gamma(x_0)=J(x_0)=0$, it
clearly follows
$$\sup_{\lambda\in ]\inf_X J, \sup_X J[}\inf_{x\in X}
(\gamma(x)-\mu\varphi(J(x)-\lambda))\leq\max\{-\epsilon,-\mu\nu\}<0$$
and so the conclusion follows in view of $(3)$.
\hfill $\bigtriangleup$
\medskip
REMARK 1. - It is clear that if a $\varphi:]-\hbox {\rm osc}_XJ,\hbox {\rm osc}_XJ[\to [0,+\infty[$ 
satisfies $(1)$ and is convex, then
it is continuous and satisfies $(2)$ too.\par
\medskip
A joint application of Theorem 1 and Proposition 1 gives\par
\medskip
THEOREM 2. -
{\it Let $X$ be a separable and reflexive real
Banach space and let $\eta, J:X\to {\bf R}$ be two $C^1$ functionals
with compact derivative and $\eta(0)=J(0)=0$. Assume also that
$J$ is bounded and non-constant, and that $\eta$ is bounded above.\par
Then,
for every
 sequentially weakly
lower semicontinuous and coercive $C^1$ functional $\psi:X\to {\bf R}$
 whose
derivative admits a continuous inverse on $X^*$ and with $\psi(0)=0$,
for every convex $C^1$ function $\varphi:
]-\hbox {\rm osc}_XJ,\hbox {\rm osc}_XJ[\to [0,+\infty[$, with
$\varphi^{-1}(0)=\{0\}$,
for which the number
$$\hat\theta=
\inf_{x\in J^{-1}({\bf R}\setminus \{0\})}
{{\psi(x)-\eta(x)}\over {\varphi(J(x))}}$$
is non-negative, and for every $\mu>\hat\theta$
there exist an
open interval $A\subseteq ]\inf_XJ,\sup_XJ[$ and a number $\rho>0$ such that,
for each $\lambda\in A$, the equation
$$\psi'(x)=\mu\varphi'(J(x)-\lambda)J'(x)+\eta'(x)$$
has at least three solutions whose norms are less than $\rho$.}\par
\smallskip
PROOF. We apply Theorem A taking $I=]\inf_XJ,\sup_XJ[$ and
$$\Psi(x,\lambda)=\psi(x)-\eta(x)-\mu\varphi(J(x)-\lambda)$$
for all $(x,\lambda)\in X\times I$.\par
Clearly, $\Psi$ is $C^1$ in $X$, continuous in $X\times I$ and concave
in $I$.
By Corollary 41.9 of [4], the functionals $\eta, J$ are sequentially weakly
continuous. Hence, for each $\lambda\in I$, the functional $\Psi(\cdot,\lambda)$
is sequentially weakly lower semicontinuous. Moreover, it is coercive, since
$\psi$ is so and $\sup_{x\in X}\max\{|J(x)|, \eta(x)\}<+\infty$.
Moreover, it is clear that, for each $\lambda\in I$, 
the derivative of the functional $\eta(\cdot)+\varphi(J(\cdot)-\lambda)$ is compact (due
to the assumptions on $\eta$ and $J$ and to the fact
that $\varphi'$ is bounded on the compact interval $[\inf_XJ,\sup_XJ]-\lambda$), and so, 
by Example 38.25 of [4], the functional
$\Psi(\cdot,\lambda)$ satisfies the Palais-Smale condition. Now, to 
realize that
condition $(a_3)$ is satisfied, we use Remark 1 and Proposition 1
with $\gamma=\psi-\eta$, observing
that $\hat\theta=\theta$ since the range of $J$ is an interval.
Then, we see that all the assumptions of Theorem
A are statisfied, and the conclusion follows in view of the chain rule.
\hfill $\bigtriangleup$
\medskip
It is worth noticing the following consequence of Theorem 2:\par
\medskip
THEOREM 3. -  
{\it Let $X$ be a separable and reflexive real
Banach space, let $J:X\to {\bf R}$ be a non-constant bounded $C^1$ functional
with compact derivative and $J(0)=0$, and let
$\psi:X\to {\bf R}$ be a
 sequentially weakly
lower semicontinuous and coercive $C^1$ functional 
 whose
derivative admits a continuous inverse on $X^*$ and with $\psi(0)=0$.
Assume that there exists $\mu>0$ such that
$$\inf_{x\in X}(\psi(x)-\mu(e^{J(x)}-1))<0\leq\inf_{x\in X}(\psi(x)-\mu J(x))
\ .\eqno{(6)}$$
Then, there exist an
open interval $A\subseteq ]\mu e^{-\sup_XJ},\mu e^{-\inf_XJ}[$ and a number
$\rho>0$ such that,
for each $\lambda\in A$, the equation
$$\psi'(x)=\lambda e^{J(x)}J'(x)$$
has at least three solutions whose norms are less than $\rho$.}\par
\smallskip
PROOF. From $(6)$, it clearly follows that
$$0\leq \inf_{x\in J^{-1}({\bf R}\setminus \{0\})}
{{\psi(x)-\mu J(x)}\over {e^{J(x)}-J(x)-1}}<\mu\ .$$
Consequently, we can apply Theorem 2 with $\eta=\mu J$ and
$\varphi(t)=e^t-t-1$, so that $\mu>\hat\theta$. Then,
there exist an open interval $B\subseteq ]\inf_XJ,\sup_XJ[$ and
a number $\rho$ such that, for each $\nu\in B$ the equation
$$\psi'(x)=\mu(e^{J(x)-\nu}-1)J'(x)+\mu J'(x)=\mu e^{-\nu}e^{J(x)}J'(x)$$
has at least three solutions whose norms are less that $\rho$. Therefore,
the conclusion follows taking
$$A=\{\mu e^{-\nu} : \nu\in B\}\ ,$$
and the proof is complete.\hfill $\bigtriangleup$\par
\medskip
{\it Proof of Theorem 1}. Let us apply Theorem 2 taking
$$X=H^1_0(\Omega)\ ,$$
$$J=J_f\ ,$$
$$\eta=J_g\ ,$$
$$\varphi=H$$
and
$$\psi(u)={{1}\over {2}}K(\|u\|^2)$$
for all $u\in X$.\par
Since $f, g\in {\cal A}$, the functionals $J_f, J_g$ are $C^1$, with compact derivative.
Since $K$ is $C^1$, increasing and coercive, the functional $\psi$ is sequentially
weakly lower semicontinuous, $C^1$ and coercive. Let us show that $\psi'$ has a continuous
inverse on $X^*$ (identified to $X$, since $X$ is a real Hilbert space).
To this end, note that the continuous function $t\to tk(t^2)$ is increasing in
$[0,+\infty[$ and onto $[0,+\infty[$. Denote by $\sigma$ its inverse and consider
the operator $T:X\to X$  defined by
$$T(v)=\cases {{{\sigma(\|v\|)}\over {\|v\|}}v & if $v\neq 0$\cr & \cr
0 & if $v=0$\ .\cr}$$
 Since $\sigma$ is continuous
and $\sigma(0)=0$, the operator $T$ is continuous in $X$. For each $u\in X\setminus \{0\}$, since $k(\|u\|^2)>0$,
we have
$$T(\psi'(u))=T(k(\|u\|^2)u)={{\sigma(k(\|u\|^2)\|u\|)}\over {k(\|u\|^2)\|u\|}}k(\|u\|^2)u={{\|u\|}\over {k(\|u\|^2)\|u\|}}k(\|u\|^2)u=u\ ,$$
as desired. Clearly, the assumptions on $h$ imply that $\varphi$ is non-negative, convex, with $\varphi^{-1}(0)=\{0\}$.
So, all the assumptions of Theorem 2 are satisfied,
and the conclusion follows.\hfill $\bigtriangleup$
\medskip
We conclude pointing out the following sample of application of Theorem 1 which is made possible by the fact
that $h$ is assumed to have the required properties on $]-\omega_f,\omega_f[$ only.\par
\medskip
EXAMPLE 1. - {\it Let $f:{\bf R}\to {\bf R}$ be a non-zero function
belonging to ${\cal A}$, with $\sup_{\bf R}|F|<+\infty$ and
let $k:[0,+\infty[\to {\bf R}$ be a continuous and non-decreasing function,
with $k(t)>0$ for all $t>0$.\par
Then, for each $\mu$ large enough,
there exist an open interval $A\subseteq
]\hbox {\rm meas}(\Omega)\inf_{\bf R}F,$\par
\noindent
$\hbox {\rm meas}(\Omega)\sup_{\bf R}F[$
and a number $\rho>0$ such that, for every $\lambda\in A$,
the problem
$$\cases {-k\left ( \int_{\Omega}|\nabla u(x)|^2dx\right )\Delta u=
\mu{{\int_{\Omega}F(u(x))dx-\lambda}\over {(\hbox {\rm meas}(\Omega) \hbox {\rm osc}_{\bf R}F)^2- 
\left ( \int_{\Omega}F(u(x))dx-\lambda\right ) ^2}} f(u) &
in $\Omega$\cr & \cr u=0 & on $\partial\Omega$\cr}$$
has at least three weak solutions whose
norms in $H^1_0(\Omega)$ are less than $\rho$.}\par
\vfill\eject
\centerline {\bf References}\par
\bigskip
\bigskip
\noindent
[1]\hskip 5pt X. L. FAN, {\it On nonlocal $p(x)$-Laplacian
Dirichlet problems}, Nonlinear Anal., {\bf 72} (2010), 3314-3323.\par
\smallskip
\noindent
[2]\hskip 5pt B. RICCERI, {\it On a three critical points theorem},
 Arch. Math. (Basel), {\bf 75} (2000), 220-226.\par
\smallskip
\noindent
[3]\hskip 5pt B. RICCERI, {\it Nonlinear eigenvalue problems}, in
{\it Handbook of Nonconvex Analysis and Applications}, D. Y. Gao and D. Motreanu eds.,
International Press, to appear.\par
\smallskip
\noindent
[4]\hskip 5pt E. ZEIDLER, {\it Nonlinear functional analysis and its
applications}, vol. III, Springer-Verlag, 1985.\par

\bigskip
\bigskip
\bigskip
\bigskip
Department of Mathematics\par
University of Catania\par
Viale A. Doria 6\par
95125 Catania\par
Italy\par
{\it e-mail address}: ricceri@dmi.unict.it

\bye